\documentclass[american,a4paper]{amsart}
\pdfoutput=1
\usepackage[T1]{fontenc}
\usepackage[utf8]{inputenc}
\usepackage{lmodern}
\usepackage{enumitem}
\usepackage{mathtools}
\usepackage{amssymb}
\usepackage{bm}
\usepackage[unicode=true,pdfusetitle,
 bookmarks=true,bookmarksnumbered=false,bookmarksopen=false,
 breaklinks=false,pdfborder={0 0 1},backref=false,colorlinks=false]{hyperref}
\usepackage{breakurl}
\usepackage[foot]{amsaddr}
\usepackage{cite}

\title[On ``Existence and uniqueness of coexistence states for an elliptic system coupled in the linear part'']{A note on ``Existence and uniqueness of coexistence states for an elliptic system coupled in the linear part'', by Hei Li-jun, Nonlinear Anal. Real World Appl. 5, 2004}
\author{L\'{e}o Girardin}
\address[L. G.]{CNRS, Institut Camille Jordan, Universit\'{e} Claude Bernard Lyon-1, 43 boulevard du 11 novembre 1918, 69622 Villeurbanne Cedex, France}
\email{leo.girardin@math.cnrs.fr}

\begin{document}
\begin{abstract}
    In this short paper I report on a paper published in Nonlinear Analysis: Real World Applications in 2004. 
    There is a major mistake early in that paper which makes most of its claims false. The class of reaction--diffusion
    systems considered in the paper has been the object of a renewed investigation in the past few years, by myself and others,
    and recent discoveries provide explicit counter-examples.
\end{abstract}

\keywords{Non-cooperative systems, KPP systems, coexistence state.}
\subjclass[2020]{35J55, 35B35.}
\maketitle

\section{Introduction}

The paper ``Existence and uniqueness of coexistence states for an elliptic system coupled in the linear part'' by Hei Li-jun 
\cite{Li-jun_2004} studies the following class of elliptic systems:
\begin{equation}\label{sys:KPP_elliptic}
    \begin{cases}
	-d_i\Delta u_i=\sum_{j=1}^n a_{ij}(x)u_j - b_i(\mathbf{u}) & \text{in }\Omega, \\
	\partial u_i/\partial \nu=0 & \text{on }\partial\Omega,
    \end{cases}
    \quad i\in\{1,2,\dots,n\},
\end{equation}
where $\Omega$ is a nonempty bounded smooth connected open set in some Euclidean space, $\partial/\partial\nu$ denotes the derivative
in the direction of the outward normal at some point on $\partial\Omega$, $n\in\mathbb{N}$, $d_1,d_2,\dots,d_n>0$, the matrix 
$\mathbf{A}=(a_{ij})_{1\leq i,j\leq n}$ is essentially positive (namely, with positive off-diagonal entries), 
$\mathbf{u}=(u_1,\dots,u_n)^{\textup{T}}$, and the $\mathcal{C}^1$ vector field $\mathbf{b}=(b_1,\dots,b_n)^{\textup{T}}$ satisfies the
following assumptions, for each $1\leq i\leq n$:
\begin{enumerate}
    \item $b_i(\mathbf{u})$ is nondecreasing, $b_i(\mathbf{0})=0$, $\partial b_i/\partial u_j(\mathbf{0})=0$, 
	for each $1\leq j\leq n$;
    \item $b_i(\mathbf{u})/u_i$ is increasing and $\lim_{u_i\to+\infty}b_i(\mathbf{u})/u_i=+\infty$ provided $u_j>0$ for each 
	$1\leq j\leq n$.
\end{enumerate}
As explained in \cite[Section 5]{Li-jun_2004}, a typical example of vector field $\mathbf{b}$ is the Lotka--Volterra competition term:
$b_i(\mathbf{u})=u_i\sum_{j=1}^n c_{ij} u_j$ with the matrix $\mathbf{C}=(c_{ij})_{1\leq i,j\leq n}$ positive entry-wise.

In order to fix the ideas, if $n=2$, then the system with Lotka--Volterra-type $\mathbf{b}$ has the following more familiar form:
\begin{equation*}
    \begin{cases}
	-d_1\Delta u_1=a_{11}(x)u_1+a_{12}(x)u_2-u_1(c_{11}u_1+c_{12}u_2), \\
	-d_2\Delta u_2=a_{21}(x)u_1+a_{22}(x)u_2-u_2(c_{21}u_1+c_{22}u_2).
    \end{cases}
\end{equation*}
Such systems are referred to as \textit{mutation--competition--diffusion systems} or \textit{stage-structured systems},
depending on the underlying biological model. They have been studied a lot in the past decades, with a renewed interest recently; 
refer for instance to
\cite{Girardin_2016_2,Morris_Borger_Crooks,Cantrell_Cosner_Yu_2018,Dockery_1998,Griette_Raoul,Bouguima_Feikh_Hennaoui_2008} 
and references therein.

After a brief reminder on linear cooperative systems in Section 2, Section 3 of \cite{Li-jun_2004} studies the nonlinear
system \eqref{sys:KPP_elliptic} and makes the following claims.
\begin{enumerate}
    \item There exists a positive solution of \eqref{sys:KPP_elliptic} if and only if the principal eigenvalue $\lambda_1$
	associated with the linearized operator at $\mathbf{u}=\mathbf{0}$, $-\operatorname{diag}(d_i\Delta)-\mathbf{A}$ 
	with Neumann boundary conditions, is negative.
    \item If $\mathbf{A}$ is symmetric ($\mathbf{A}=\mathbf{A}^{\textup{T}}$) and $\lambda_1<0$, then the positive solution is unique.
\end{enumerate}
Stability considerations follow in Section 4.

Although the first claim is consistent with other works on this class of systems, the second one is contradictory. It turns out
that even the proof of existence when $\lambda_1<0$ is false, as written. I will explain the error and what can be salvaged
in the following sections of the present note.

\section{The uniqueness result cannot possibly be true}

In \cite{Girardin_2017}, I gave the following explicit counter-example showing that positive coexistence states are not, 
in general, unique, even with a symmetry assumption on $\mathbf{A}$:
\begin{equation}\label{sys:counterexample}
    \begin{cases}
	-d_1\Delta u_1=(1-1/5)u_1+(1/5)u_2-(1/10)u_1(u_1+9u_2), \\
	-d_2\Delta u_2=(1/5)u_1+(1-1/5)u_2-(1/10)u_2(9u_1+u_2).
    \end{cases}
\end{equation}
Indeed, it can be directly verified that the constant positive solutions are in this case:
\[
    \begin{pmatrix}1\\1\end{pmatrix},\quad\begin{pmatrix}3-\sqrt{15/2}\\3+\sqrt{15/2}\end{pmatrix},\quad\begin{pmatrix}3+\sqrt{15/2}\\3-\sqrt{15/2}\end{pmatrix}.
\]

More recently, Cantrell, Cosner and Yu \cite{Cantrell_Cosner_Yu_2018} actually showed that, in the two-species case $n=2$,
any system of the form \eqref{sys:KPP_elliptic} with Lotka--Volterra competition which is bistable competitive in the absence 
of mutations ($a_{12}=a_{21}=0$) will remain bistable competitive
with sufficiently small mutations, in the sense that the two stable semi-trivial steady states of the Lotka--Volterra competitive 
system are displaced by small positive off-diagonal terms $a_{12}, a_{21}$ inside the positive cone $\{u_1\geq 0, u_2\geq 0\}$ and
not outside of it. The preceding explicit counter-example \eqref{sys:counterexample} is a particular case of this general result.

The existence result depending on the sign of $\lambda_1$ is, on the contrary, consistent with the literature (for instance,
\cite[Theorem 1.4]{Girardin_2016_2}). 

\section{In fact, all proofs of the paper are false}

Nevertheless, after a thorough inspection of the proofs in \cite{Li-jun_2004}, it appears that even the proof of existence
is false. In fact, the main idea of the paper, namely that \eqref{sys:KPP_elliptic} is actually cooperative, is wrong, and since
all proofs rely at least partially on this idea (monotonicity arguments, comparison principles and stability theory for monotone 
systems), all results are false.

The system \eqref{sys:KPP_elliptic} is indeed cooperative at the origin $\mathbf{u}=\mathbf{0}$ but the preservation of this 
property away from the origin strongly depends on the specific choice of $\mathbf{b}$ and, although it can be true 
\cite{Cantrell_Cosner_Yu_2018}, it is in general false \cite{Cantrell_Cosner_Yu_2018,Girardin_2018}.

More precisely, the error leading to this false nonlinear comparison principle is in the proof of Theorem 3.2 of \cite{Li-jun_2004}.
At the very beginning of the proof, the author states the following claim: fix $\underline{\mathbf{u}}\gg\mathbf{0}$,
$\overline{\mathbf{u}}\gg\underline{\mathbf{u}}$\footnote{As stated in \cite{Li-jun_2004}, the claim only assumes 
$\overline{\mathbf{u}}\geq\underline{\mathbf{u}}\geq\mathbf{0}$, however the actual super- and sub-solutions constructed
later on in \cite{Li-jun_2004} all satisfy $\overline{\mathbf{u}}\gg\underline{\mathbf{u}}\gg\mathbf{0}$. I claim, without
proof for the sake of brevity, that the construction of a counter-example when the inequalities are large instead of strict 
is cumbersome but still possible.}, and let $\mathbf{u}_1$, $\mathbf{u}_2$ be two arbitrary vectors
such that $\underline{\mathbf{u}}\leq \mathbf{u}_2\leq \mathbf{u}_1 \leq \overline{\mathbf{u}}$. Then,
provided the real number $M>0$ is sufficiently large in a way that does not depend on $\mathbf{u}_1$ and $\mathbf{u}_2$,
\[
    M\mathbf{u}_1-\mathbf{b}(\mathbf{u}_1)-M\mathbf{u}_2+\mathbf{b}(\mathbf{u}_2)\geq 0.
\]
Although this is obviously true in the scalar setting $n=1$ (with $M$ the Lipschitz constant of $\mathbf{b}$), this is false
in higher dimensions $n\geq 2$, as can be easily verified with the following counter-example.
Let $\mathbf{v}=\underline{\mathbf{u}}$, $\mathbf{u}_2=\mathbf{v}$,
$\mathbf{u}_1=\mathbf{v}+\varepsilon\mathbf{e}_1=\mathbf{v}+(\varepsilon,0,\dots,0)^{\textup{T}}$ with 
$\varepsilon>0$ so small that $\mathbf{u}_1\leq\overline{\mathbf{u}}$. Then for each $i\geq 2$ and any $M>0$,
\[
    \left[M\mathbf{u}_1-\mathbf{b}(\mathbf{u}_1)-M\mathbf{u}_2+\mathbf{b}(\mathbf{u}_2)\right]_i
    = b_i(\mathbf{v})-b_i(\mathbf{v}+\varepsilon\mathbf{e}_1)
    = v_i\left(\frac{b_i(\mathbf{v})}{v_i}-\frac{b_i(\mathbf{v}+\varepsilon\mathbf{e}_1)}{v_i}\right)<0.
\]
(Even though the assumption of strict monotonicity of $b_i(\mathbf{u})/u_i$ is ambiguously stated in \cite{Li-jun_2004}, 
the previous strict inequality is clear in the aforementioned Lotka--Volterra competition case that the author clearly had in mind.)
Everything that follows in the paper depends on this false claim.

\section{What can be salvaged}

Theorem 3.4 in \cite{Li-jun_2004} states the uniqueness of the coexistence state under a symmetry assumption on the matrix $\mathbf{A}$.
Although it is false in general as explained before, the variational argument there does not directly rely upon the comparison
principle and can be used to prove that two distinct coexistence states $\mathbf{u}$ and $\mathbf{v}$ cannot be compared 
(namely, both $\mathbf{u}\geq \mathbf{v}$ and $\mathbf{u}\leq \mathbf{v}$ are false). This property
is indeed satisfied by the above counter-example \eqref{sys:counterexample}. 

\bibliographystyle{plain}
\bibliography{ref}

\end{document}